\documentclass[12pt]{article}

\newcommand{\qed}{{\hfill\rule{4pt}{7pt}}}

\makeatletter
\newfont{\footsc}{cmcsc10 at 8truept}
\newfont{\footbf}{cmbx10 at 8truept}
\newfont{\footrm}{cmr10 at 10truept}
\makeatother
\usepackage{amsmath}
\usepackage{amssymb}
\usepackage{amsfonts}
\usepackage{latexsym}

\newtheorem{thm}{Theorem}[section]

\newtheorem{cor}[thm]{Corollary}

\newtheorem{lem}[thm]{Lemma}

\def\pf{\noindent {\it Proof.} }
\numberwithin{equation}{section}

\title{A $q$-Analogue of Faulhaber's Formula for Sums of Powers}

\author{Victor J. W. Guo and Jiang Zeng\\
\small Institut Camille Jordan, Universit\'e Claude Bernard (Lyon I)\\[-0.8ex]
\small F-69622 Villeurbanne Cedex, France\\[-0.8ex]
\small \texttt{\tt jwguo@eyou.com, zeng@math.univ-lyon1.fr}}

\date{\small
Submitted: Jan 25, 2005; Accepted: Aug 16, 2005; Published: ******, 2005}

\begin{document}
\maketitle


\begin{abstract}
Let
 $$
 S_{m,n}(q):=\sum_{k=1}^{n}\frac{1-q^{2k}}{1-q^2}
 \left(\frac{1-q^k}{1-q}\right)^{m-1}q^{\frac{m+1}{2}(n-k)}.
 $$
Generalizing the formulas of Warnaar and Schlosser, we prove that
there exist polynomials $P_{m,k}(q)\in\mathbb{Z}[q]$ such that
 $$
 S_{2m+1,n}(q) =\sum_{k=0}^{m}(-1)^kP_{m,k}(q)
 \frac{(1-q^n)^{m+1-k}(1-q^{n+1})^{m+1-k}q^{kn}}
 {(1-q^2)(1-q)^{2m-3k}\prod_{i=0}^{k}(1-q^{m+1-i})},
$$
and solve a problem raised by
Schlosser. We also show that there is a similar
formula for the following $q$-analogue of alternating sums of powers:
 $$
 T_{m,n}(q):=\sum_{k=1}^{n}(-1)^{n-k}
 \left(\frac{1-q^k}{1-q}\right)^{m}q^{\frac{m}{2}(n-k)}.
 $$
\end{abstract}

\section{Introduction}
In the early 17th century Faulhaber~\cite{Faulhaber} computed the
sums of powers $1^m+2^m+\cdots +n^m$ up to $m=17$ and realized that for
odd $m$, it is not just a polynomial in $n$ but a polynomial in the
triangular number $N=n(n+1)/2$.
A good account of Faulhaber's work was given by Knuth~\cite{Knuth}.
For example, for $m=1, \ldots, 5$, Faulhaber's formulas
read as follows:
\begin{align*}
1^1+2^1+\cdots+n^1 &=N,\qquad N=(n^2+n)/2\,;\\[5pt]
1^2+2^2+\cdots +n^2&= \frac{2n+1}{3}N\,;
\end{align*}
\begin{align*}
1^3+2^3+\cdots+n^3 &=N^2;\\[5pt]
1^4+2^4+\cdots +n^4 &=\frac{2n+1}{5}(2N^2-\frac{1}{3}N)\,;\\[5pt]
1^5+2^5+\cdots +n^5 &=\frac{1}{3}(4N^3-N^2).
\end{align*}
Recently, the problem of $q$-analogues of the sums of powers has
attracted the attention of several authors
\cite{GH,Warnaar,Schlosser}, who found, in particular,
 $q$-analogues of the Faulhaber formula corresponding to
$m=1,2,\ldots, 5$.
More precisely, setting
\begin{align}
S_{m,n}(q)=\sum_{k=1}^{n}\frac{1-q^{2k}}{1-q^2}
\left(\frac{1-q^k}{1-q}\right)^{m-1}q^{\frac{m+1}{2}(n-k)},
\label{eq:qpower}
\end{align}
Warnaar~\cite{Warnaar} (for $m=3$) and  Schlosser~\cite{Schlosser} found
 the following formulas for the  $q$-analogues of the sums of consecutive integers, squares,
cubes, quarts and quints:
\begin{align}
S_{1,n}(q)&=\frac{(1-q^n)(1-q^{n+1})}{(1-q)(1-q^2)},\label{eq:f1}\\[5pt]
S_{2,n}(q)&=\frac{(1-q^n)(1-q^{n+1})(1-q^{n+\frac{1}{2}})}{(1-q)(1-q^2)(1-q^{\frac32})},
\\[5pt]
S_{3,n}(q)&=\frac{(1-q^n)^2(1-q^{n+1})^2}{(1-q)^2(1-q^2)^2},\\[5pt]
S_{4,n}(q)&=\frac{(1-q^n)(1-q^{n+1})(1-q^{n+\frac{1}{2}})}{(1-q)(1-q^2)(1-q^{\frac{5}{2}})}
\Bigg[\frac{(1-q^n)(1-q^{n+1})}{(1-q)^2}-\frac{1-q^{\frac{1}{2}}}{1-q^{\frac{3}{2}}}
q^{n}\Bigg], \label{eq:f4}\\[5pt]
S_{5,n}(q)&=\frac{(1-q^n)^2(1-q^{n+1})^2}{(1-q)^2(1-q^2)(1-q^3)}
\Bigg[\frac{(1-q^n)(1-q^{n+1})}{(1-q)^2}-\frac{1-q}{1-q^2}q^{n}\Bigg]. \label{eq:f5}
\end{align}
Notice that the above formulas have the same pattern that each
summand on the right-hand side has no pole at $q=1$, and
so reduce directly to Faulhaber's corresponding formulas when
$q\to 1$.

At the end of his paper, Schlosser~\cite{Schlosser} speculated on the existence of
a general formula for $S_{m,n}(q)$, and left it as an open problem.
It is the purpose of this paper to provide such a general formula,
which turns out to be a $q$-analogue of the Faulhaber formula for the sums of powers.
More precisely, we prove the following results:
\begin{thm}\label{thm:pmn}
For $m,n\in \mathbb{N}$, there exist polynomials $P_{m,k}(q)\in\mathbb{Z}[q]$ such that
\begin{align*}
S_{2m+1,n}(q) &=\sum_{k=0}^{m}(-1)^kP_{m,k}(q)
\frac{(1-q^n)^{m+1-k}(1-q^{n+1})^{m+1-k}q^{kn}}
{(1-q^2)(1-q)^{2m-3k}\prod_{i=0}^{k}(1-q^{m+1-i})},
\end{align*}
where $P_{m,s}(q)$ are the $q$-Faulhaber coefficients given by
\begin{align}
P_{m,s}(q)&=\frac{\prod_{j=0}^{s}(1-q^{m+1-j})}{(1-q)^{3s}}\sum_{k=0}^{s}
\frac{(-1)^{s-k}}{1-q^{m+1-k}}
\left[{2m\choose k }-{2m\choose k -2}\right] \nonumber\\[5pt]
&\quad{}\times\sum_{i=0}^{s-k}\frac{m-k+1}{m-s+1} {m-s+i\choose i}{m-k-i\choose s-k-i}
q^{s-k-i}.
\label{eq:cor1}
\end{align}
\end{thm}

\begin{thm}\label{thm:qmn}
For $m,n\in \mathbb{N}$, there exist polynomials $Q_{m,k}(q)\in\mathbb{Z}[q]$ such that
\begin{align*}
S_{2m,n}(q)
=\sum_{k=0}^{m}(-1)^k Q_{m,k}(q^{\frac{1}{2}})
\frac{(1-q^{n+\frac{1}{2}})(1-q^n)^{m-k}(1-q^{n+1})^{m-k}(1-q^{\frac{1}{2}})^kq^{kn}}
{(1-q^2)(1-q)^{2m-2k-1}\prod_{i=0}^{k}(1-q^{m-i+\frac{1}{2}})}.
\end{align*}
Furthermore, we have
\begin{align}
Q_{m,s}(q)&=\frac{\prod_{j=0}^{s}(1-q^{2m-2j+1})}{(1-q)^{s}(1-q^2)^{2s}}
\sum_{k=0}^{s}\frac{(-1)^{s-k}}{1-q^{2m-2k+1}}\left[{2m-1\choose k}-{2m-1\choose k-2}\right]
\nonumber\\[5pt]
&\quad{}\times\sum_{i=0}^{s-k} {m-s+i\choose i}
\bigg[{m-k-i\choose s-k-i}q^{2s-2k-2i}+{m-k-i-1\choose s-k-i-1}q^{2s-2k-2i-1}\bigg].  \label{eq:cor2}
\end{align}
\end{thm}

Next we consider a $q$-analogue of the alternating sums
$T_{m,n}=\sum_{k=1}^{n}(-1)^{n-k}k^{m}$. Note that Gessel and
Viennot~\cite{GV} proved that $T_{2m,n}$ can be written as a
polynomial in $n(n+1)$ whose coefficients are the Sali\'e
coefficients and Schlosser~\cite{Schlosser} gave some $q$-analogues
of $T_{m,n}$ only for $m\leq 4$. Let
\begin{align}
T_{m,n}(q)=\sum_{k=1}^{n}(-1)^{n-k}
\left(\frac{1-q^k}{1-q}\right)^{m}q^{\frac{m}{2}(n-k)}.
\label{eq:altersum}
\end{align}
We have the following $q$-analogue of Gessel-Viennot's result for $T_{m,n}$.
\begin{thm}\label{thm:salie}
For $m,n\in \mathbb{N}$, there exist polynomials $G_{m,k}(q)\in\mathbb{Z}[q]$ such that
\begin{align}
T_{2m,n}(q) &=\sum_{k=0}^{m-1}(-1)^k G_{m,k}(q)
\frac{(1-q^n)^{m-k}(1-q^{n+1})^{m-k}q^{kn}}
{(1-q)^{2m-2k}\prod_{i=0}^{k}(1+q^{m-i})}, \label{eq:t2mnq}
\end{align}
where $G_{m,k}(q)$ are the $q$-Sali\'e coefficients given by
\begin{align*}
G_{m,s}(q)&=\frac{\prod_{j=0}^{s}(1+q^{m-j})}{(1-q)^{2s}}
             \sum_{k=0}^{s}\frac{(-1)^{s-k}}{1+q^{m-k}}{2m\choose k}\\[5pt]
&\quad{}\times \sum_{i=0}^{s-k}\frac{m-k}{m-s}
{m-s+i-1\choose i}{m-k-i-1\choose s-k-i}q^{s-k-i}.
\end{align*}
\end{thm}

\begin{thm}\label{thm:salodd}
For $m,n\in \mathbb{N}$, there exist polynomials $H_{m,k}(q)\in\mathbb{Z}[q]$ 
such that
\begin{align*}
T_{2m-1,n}(q) &=(-1)^{m+n}H_{m,m-1}(q^{\frac12})\frac{q^{mn-\frac{n}{2}}}
{(1+q^{\frac12})^{m}\prod_{i=0}^{m-1}(1+q^{m-i-\frac12})} \\[5pt]
&\quad+\frac{1-q^{n+\frac12}}{1-q^{\frac12}}\sum_{k=0}^{m-1}(-1)^k
\frac{H_{m,k}(q^{\frac12}) (1-q^n)^{m-k-1}(1-q^{n+1})^{m-k-1}q^{kn}}
{(1-q)^{2m-2k-2}(1+q^{\frac12})^{k+1}\prod_{i=0}^{k}(1+q^{m-i-\frac12})}.
\end{align*}
Furthermore, we have
\begin{align*}
H_{m,s}(q)&=\frac{\prod_{j=0}^{s}(1+q^{2m-2j-1})}{(1+q)^{s}(1-q)^{2s}}
\sum_{k=0}^{s}\frac{(-1)^{s-k}}{1+q^{2m-2k-1}}{2m-1\choose k}
\sum_{i=0}^{s-k} {m-s+i-1\choose i}\\[5pt]
&\quad{}\times
\bigg[{m-k-i-1\choose s-k-i}q^{2s-2k-2i}+{m-k-i-2\choose s-k-i-1}
q^{2s-2k-2i-1}\bigg].
\end{align*}
\end{thm}

Schlosser~\cite{Schlosser} derives his formulas
from the machinery of basic hypergeometric series. For example,
for the $q$-analogues of the sums of quarts and quints, he first
specializes Bailey's terminating very-well-poised balanced
${}_{10}\phi_{9}$ transformation \cite[Appendix (III.28)]{GR} and
then applies  the terminating very-well-poised ${}_6\phi_{5}$
\cite[Appendix (II.21)]{GR} on one side of the identity to
establish a ``master identity." In contrast to his proof, our method is
self-contained and of elementary nature.

We first establish some elementary algebraic identities in Section~\ref{preli},
and prove Theo-rems~\ref{thm:pmn}--\ref{thm:salodd} in Section~\ref{appl}.
We then apply our theorems to compute the polynomials
$P_{m,s}(q)$, $Q_{m,s}(q)$, $G_{m,s}(q)$, and $H_{m,s}(q)$ for
small $m$ in Section~\ref{compute}
and obtain summation formulas of \eqref{eq:qpower} for
$m\leq 11$. Section~\ref{sec:open} contains some further extensions of
these summation formulas.

\section{Some Preliminary Lemmas}\label{preli}
The following is our first step towards our summation formula for $S_{m,n}(q)$.
\begin{lem}\label{lem:main}
For $m,n\in\mathbb{N}$, we have
\begin{align}
S_{m,n}(q)&=\sum_{r=0}^{\lfloor \frac m2\rfloor}(-1)^r
\left[{m-1\choose r}-{m-1\choose
r-2}\right]
\frac{(1-q^{(\frac{m+1}{2}-r)n})(1+(-1)^mq^{(\frac{m+1}{2}-r)(n+1)})q^{rn}}
{(1-q^2)(1-q)^{m-1}(1-q^{\frac{m+1}{2}-r})}.\label{eq:1sum}
\end{align}
\end{lem}
\pf By definition, $(1-q^2)(1-q)^{m-1}S_{m,n}(q)$ is equal to
\begin{align}
&\hskip -3mm 
\sum_{k=1}^{n}(1-q^{2k})(1-q^k)^{m-1}q^{\frac{m+1}{2}(n-k)}\nonumber\\[5pt]
&=\sum_{k=1}^{n}(1-q^{2k})q^{\frac{m+1}{2}(n-k)}
\sum_{r=0}^{m-1}{m-1\choose r}(-1)^rq^{kr}  \nonumber
\end{align}
\begin{align}
&=\sum_{r=0}^{m-1}{m-1\choose r}(-1)^r
\sum_{k=1}^{n}(q^{\frac{m+1}{2}n+(r-\frac{m+1}{2})k}-q^{\frac{m+1}{2}n
+(r-\frac{m-3}{2})k})\nonumber\\[5pt]
&=\sum_{r=0}^{m+1}(-1)^r\left[{m-1\choose r}-{m-1\choose
r-2}\right]
\sum_{k=1}^{n}q^{\frac{m+1}{2}n+(r-\frac{m+1}{2})k}\nonumber\\[5pt]
&=\sum_{\substack{r=0\\ r\neq \frac{m+1}2}}^{m+1}(-1)^r
\left[{m-1\choose r}-{m-1\choose r-2}\right]
\frac{q^{\frac{m+1}{2}(n-1)+r}-q^{r(n+1)-\frac{m+1}{2}}}
{1-q^{r-\frac{m+1}{2}}}.\label{eq:bigsum}
\end{align}
Splitting the last summation into two parts corresponding to $r$
ranging from 0 to ${\lfloor \frac m2\rfloor}$
and from ${\lfloor \frac{m+1}{2}\rfloor+1}$ to $m+1$, respectively.
Replacing $r$ by $m+1-r$ in the second
one  we can rewrite \eqref{eq:bigsum} as follows:
\begin{align*}
&\sum_{r=0}^{\lfloor \frac m2\rfloor}(-1)^r\left[{m-1\choose r}
-{m-1\choose r-2}\right] \\[5pt]
&\quad{}\times\left[\frac{q^{\frac{m+1}{2}(n-1)+r}-q^{r(n+1)-\frac{m+1}{2}}}
{1-q^{r-\frac{m+1}{2}}}+
(-1)^m\frac{q^{\frac{m+1}{2}(n+1)-r}-q^{\frac{m+1}{2}(2n+1)-r(n+1)}}
{1-q^{\frac{m+1}{2}-r}}\right].
\end{align*}
After simplification we get \eqref{eq:1sum}. \qed

\medskip
\noindent{\it Remark.} When $m$ is even, since
$$
\sum_{r=0}^{\frac m2}(-1)^r \left[{m-1\choose r}-{m-1\choose r-2}\right]=0,
$$
we can rewrite $S_{m,n}(q)$ as
\begin{equation}
S_{m,n}(q)=\sum_{r=0}^{\frac m2}(-1)^r \left[{m-1\choose r}-
{m-1\choose r-2}\right] \frac{(1-q^{(\frac{m+1}{2}-r)(2n+1)})q^{nr}}
{(1-q^2)(1-q)^{m-1}(1-q^{\frac{m+1}{2}-r})}. \label{eq:oddq}
\end{equation}

\begin{lem}\label{lem:sum}
For $m,n\geq 1$, we have
\begin{align}
T_{2m,n}(q)=\sum_{r=0}^{m-1}(-1)^r {2m\choose r}
\frac{(1-q^{n(m-r)})(1-q^{(n+1)(m-r)})q^{rn}}
{(1-q)^{2m}(1+q^{m-r})}.\label{eq:lem-tmn}
\end{align}
\end{lem}
\pf By \eqref{eq:altersum} we have
$$
(1-q)^{2m}T_{2m,n}(q)=\sum_{k=1}^{n}(1-q^k)^{2m}q^{m(n-k)}(-1)^{n-k}.
$$
Expanding $(1-q^k)^{2m}$ by the binomial theorem and exchanging the summation order,
we obtain
\begin{align}
(1-q)^{2m}T_{2m,n}(q)
=\sum_{r=0}^{2m}(-1)^r{2m\choose r}q^{rn}\frac{1-(-q^{m-r})^n}{1+q^{m-r}}.
\label{eq:t2mn0}
\end{align}
Substituting $r$ by $2m-r$ on the right-hand side of \eqref{eq:t2mn0} yields
\begin{align}
(1-q)^{2m}T_{2m,n}(q)&=\frac{1}{2}\sum_{r=0}^{2m}(-1)^r{2m\choose r}
\left[q^{rn}\frac{1-(-q^{m-r})^n}{1+q^{m-r}}
+q^{(2m-r)n}\frac{1-(-q^{r-m})^n}{1+q^{r-m}}\right] \nonumber\\[5pt]
&=\frac{1}{2}\sum_{r=0}^{2m}(-1)^r{2m\choose r}q^{rn}
\frac{(1-(-1)^nq^{n(m-r)})(1-(-1)^nq^{(n+1)(m-r)})}{1+q^{m-r}} \nonumber\\[5pt]
&=\frac{1}{2}\sum_{r=0}^{2m}(-1)^r{2m\choose r}q^{rn}
\frac{(1-q^{n(m-r)})(1-q^{(n+1)(m-r)})}{1+q^{m-r}}. \label{eq:sym-1n}
\end{align}
The last equality holds because
$$
\sum_{r=0}^{2m}(-1)^r{2m\choose r}q^{rn}\frac{q^{n(m-r)}+q^{(n+1)(m-r)}}{1+q^{m-r}}=0.
$$
Splitting  the sum
in  \eqref{eq:sym-1n} as $\sum_{r=0}^{m-1}+\sum_{r=m+1}^{2m}$
and substituting $r$ by $2m-r$ in the second sum, we complete the proof. \qed

\medskip

Similarly, we can show that $2(1-q)^{2m-1}T_{2m-1,n}(q)$ is equal to
\begin{align*}
&\hskip -3mm \sum_{r=0}^{2m-1}(-1)^r{2m-1\choose r}q^{rn}
\frac{(1-(-1)^nq^{n(m-r-\frac12)})(1+(-1)^nq^{(n+1)(m-r-\frac12)})}
{1+q^{m-r-\frac12}} \nonumber\\[5pt]
&=\sum_{r=0}^{2m-1}(-1)^{r}{2m-1\choose r}\left[(-1)^{n+1}q^{(m-\frac12)n}
\frac{1-q^{m-r-\frac12}}{1+q^{m-r-\frac12}}
+q^{rn}\frac{1-q^{(2n+1)(m-r-\frac12)}}{1+q^{m-r-\frac12}}\right].
\end{align*}
This establishes immediately the following lemma:
\begin{lem}\label{lem:todd}
For $m,n\geq 1$, we have
\begin{align}
T_{2m-1,n}(q)&=\sum_{r=0}^{m-1}(-1)^{n+r+1}{2m-1\choose r}
\frac{(1-q^{m-r-\frac12})q^{(m-\frac12)n}}{(1-q)^{2m-1}(1+q^{m-r-\frac12})} \nonumber\\[5pt]
&\quad+\sum_{r=0}^{m-1}(-1)^r{2m-1\choose r}
\frac{(1-q^{(2n+1)(m-r-\frac12)})q^{rn}}{(1-q)^{2m-1}(1+q^{m-r-\frac12})}.
\label{eq:lem-todd}
\end{align}
\end{lem}

The second ingredient of our approach is the following identity, of which we
shall give two proofs.

\begin{thm}\label{thm:mrsxy}
For $m\in \mathbb{N}$, we have
\begin{align}
\frac{1-x^{m+1}y^{m+1}}{(1-xy)(1-x)^m(1-y)^m}
=\sum_{r=0}^{m}\sum_{s=o}^{m-r}{m-r\choose s}{m-s\choose r}
\frac{x^ry^s}{(1-x)^{r+s}(1-y)^{r+s}}. \label{eq:mrsxy}
\end{align}
\end{thm}

\noindent{\it First Proof.}  Replacing $s$ by $m-r-s$, the right-hand side of
\eqref{eq:mrsxy} may be written as
\begin{align}
\sum_{r=0}^{m}\sum_{s=o}^{m-r}{m-r\choose s}{r+s\choose r}
\frac{x^ry^{m-r-s}}{(1-x)^{m-s}(1-y)^{m-s}}. \label{eq:mrsbino}
\end{align}
Consider the generating function of \eqref{eq:mrsbino}. We have
\begin{align}
&\hskip -3mm 
\sum_{m=0}^\infty\sum_{r=0}^{m}\sum_{s=o}^{m-r}{m-r\choose s}{r+s\choose r}
\frac{x^ry^{m-r-s}}{(1-x)^{m-s}(1-y)^{m-s}} t^m \nonumber\\[5pt]
&=\sum_{r=0}^{\infty}\sum_{s=0}^{\infty}{r+s\choose r} 
x^r\sum_{m=r+s}^{\infty}{m-r\choose s}
\frac{y^{m-r-s}}{(1-x)^{m-s}(1-y)^{m-s}} t^m
\nonumber\\[5pt]
&=\sum_{r=0}^{\infty}\sum_{s=0}^{\infty}{r+s\choose r} 
\frac{x^r t^{r+
s}}{(1-x)^{r}(1-y)^{r}}
\left(1-\frac{yt}{(1-x)(1-y)}\right)^{-s-1}
\nonumber\\[5pt]
&=\frac{(1-x)(1-y)}{(1-x)(1-y)-yt}
\left(1-\frac{xt}{(1-x)(1-y)}-\frac{t(1-x)(1-y)}{(1-x)(1-y)-yt}\right)^{-1}
\nonumber\\[5pt]
&=\frac{(1-x)^2(1-y)^2}{[(1-x)(1-y)-xyt][(1-x)(1-y)-t]}, \nonumber
\end{align}
which is equal to the generating function of the left-hand side of
\eqref{eq:mrsxy}. \qed

\medskip

\noindent{\it Second Proof.}
Let
\[
\begin{cases}
x=u(1-x)(1-y),\\
y=v(1-x)(1-y).
\end{cases}
\]
We want to expand
\[
f(x,y)=\frac{1-x^{m+1}y^{m+1}}{(1-xy)(1-x)^m(1-y)^m}
\]
as a series in $u$ and $v$. By Lagrange's inversion formula (see, for example, \cite[p.~21]{GJ}),
\[
f(x,y)=\sum_{r,s\geq 0}u^r v^s [x^ry^s]
\left\{\frac{1-x^{m+1}y^{m+1}}{(1-xy)(1-x)^{m-r-s}(1-y)^{m-r-s}}\Delta\right\},
\]
where $[x^ry^s]F(x,y)$ denotes the coefficient of $x^ry^s$
in the power series $F(x,y)$, and where $\Delta$ is the determinant given by
\[
\Delta=\left|
\begin{array}{rr}
\displaystyle 1+\frac{x}{1-x}& \displaystyle \frac{y}{1-y} \\[10pt]
\displaystyle \frac{x}{1-x}  & \displaystyle 1+\frac{y}{1-y}
\end{array}\right|
=\frac{1-xy}{(1-x)(1-y)}.
\]
So,
\begin{equation}\label{eq:key}
f(x,y)=\sum_{r,s\geq 0}u^r v^s [x^ry^s]
\left\{\frac{1-x^{m+1}y^{m+1}}{(1-x)^{m-r-s+1}(1-y)^{m-r-s+1}}\right\}.
\end{equation}
Since
\begin{align*}
&(1-z)^{-\alpha}=\sum_{k=0}^{\infty}{\alpha+k-1\choose k}z^k,
\end{align*}
we have
\[
[x^r y^s]\left\{(1-x)^{-(m-r-s+1)}(1-y)^{-(m-r-s+1)}\right\}
={m-s\choose r}{m-r\choose s},
\]
and
\begin{align*}
&\hskip -3mm
[x^r y^s]\left\{(1-x)^{-(m-r-s+1)}(1-y)^{-(m-r-s+1)}x^{m+1}y^{m+1}\right\}\\[5pt]
&=\begin{cases}
0,&\text{if $r\leq m$ or $s\leq m$,} \\
(-1)^{r+s}{r+s-m-1\choose s}{r+s-m-1\choose r}, &\text{if $r,s\geq m+1$.}
\end{cases}
\end{align*}
But, it is easy to see that
\[
{m-s\choose r}{m-r\choose s}=(-1)^{r+s}{r+s-m-1\choose s}{r+s-m-1\choose r}.
\]
Substituting these into \eqref{eq:key} yields
\begin{align*}
f(x,y)=\sum_{0\leq r,s\leq m}u^r v^s {m-s\choose r}{m-r\choose s}.\tag*{\qed}
\end{align*}

\begin{cor} For $m\in \mathbb{N}$, we have
\begin{align}
&\hskip -3mm\sum_{r=0}^{m}\sum_{s=o}^{m-r}{m-r+1\choose s}{m-s\choose r}
\frac{x^ry^s}{(1-x)^{r+s}(1-y)^{r+s}} \nonumber\\[5pt]
&=\frac{1-x^{m+2}y^{m+2}-x(1-x^{m+1}y^{m+1})-(1-xy)y^{m+1}}{(1-xy)(1-x)^{m+1}(1-y)^{m+1}}.
\label{eq:mrsxy0}
\end{align}
\end{cor}
\pf Replacing $m$ and $r$ by $m-1$ and $r-1$ respectively in \eqref{eq:mrsxy}, we obtain
\begin{align}
\sum_{r=1}^{m}\sum_{s=o}^{m-r}{m-r\choose s}{m-s-1\choose r-1}
\frac{x^ry^s}{(1-x)^{r+s}(1-y)^{r+s}}
=\frac{x(1-x^{m}y^{m})}{(1-xy)(1-x)^{m}(1-y)^{m}}. \label{eq:mrsxy2}
\end{align}
Combining \eqref{eq:mrsxy} and \eqref{eq:mrsxy2}, we get
\begin{align}
\sum_{r=0}^{m}\sum_{s=o}^{m-r}{m-r\choose s}{m-s-1\choose r}
\frac{x^ry^s}{(1-x)^{r+s}(1-y)^{r+s}}
=\frac{1-x^{m+1}y^{m+1}-x(1-x^{m}y^{m})}{(1-xy)(1-x)^{m}(1-y)^{m}}.
\label{eq:mrsxy3}
\end{align}
Replacing $m$ by $m+1$ in \eqref{eq:mrsxy3}, we have
\begin{align}
&\hskip -3mm\sum_{r=0}^{m+1}\sum_{s=o}^{m-r+1}{m-r+1\choose s}{m-s\choose r}
\frac{x^ry^s}{(1-x)^{r+s}(1-y)^{r+s}} \nonumber\\[5pt]
&=\frac{1-x^{m+2}y^{m+2}-x(1-x^{m+1}y^{m+1})}{(1-xy)(1-x)^{m+1}(1-y)^{m+1}}.
\label{eq:mrsxy4}
\end{align}
Note that when $r=m+1$, ${m-s\choose r}=0$, and when $s=m-r+1$,
${m-s\choose r}={r-1\choose r}$ is equal to $1$ if $r=0$ and $0$
otherwise. Moving the term
$\frac{y^{m+1}}{(1-x)^{m+1}(1-y)^{m+1}}$ of \eqref{eq:mrsxy4} from
the left-hand side to the right-hand side, we obtain
\eqref{eq:mrsxy0}. \qed

\medskip

Interchanging $r$ and $s$, and $x$ and $y$ in \eqref{eq:mrsxy0}, we get
\begin{align}
&\hskip -3mm \sum_{r=0}^{m}\sum_{s=o}^{m-r}{m-r\choose s}{m-s+1\choose r}
\frac{x^ry^s}{(1-x)^{r+s}(1-y)^{r+s}} \nonumber\\[5pt]
&=\frac{1-x^{m+2}y^{m+2}-y(1-x^{m+1}y^{m+1})-(1-xy)x^{m+1}}{(1-xy)(1-x)^{m+1}(1-y)^{m+1}}.
\label{eq:mrsxy00}
\end{align}

\begin{cor}\label{thm:1x1yik}For $m\in \mathbb{N}$, we have
\begin{align}
\frac{(1-x^{m+1})(1-y^{m+1})}{(1-x)^{m+1}(1-y)^{m+1}}
=\sum_{r=0}^{m}\sum_{s=o}^{m-r}\frac{m+1}{m+1-r-s}{m-r\choose s}{m-s\choose r}
\frac{x^ry^s}{(1-x)^{r+s}(1-y)^{r+s}}. \label{eq:1x1yrs}
\end{align}
\end{cor}

\pf Note that
\begin{align}
&\hskip -3mm \frac{m+1}{m+1-r-s}{m-r\choose s}{m-s\choose r} \nonumber\\[5pt]
&={m-r+1\choose s}{m-s\choose r}+{m-r\choose s}{m-s+1\choose r}
-{m-r\choose s}{m-s\choose r}. \label{eq:ppos}
\end{align}
Hence, from \eqref{eq:mrsxy}, \eqref{eq:mrsxy0} and \eqref{eq:mrsxy00}
it follows that
\begin{align*}
&\hskip -3mm\sum_{r=0}^{m}\sum_{s=o}^{m-r}\frac{m+1}{m+1-r-s}{m-r\choose s}{m-s\choose r}
\frac{x^ry^s}{(1-x)^{r+s}(1-y)^{r+s}} \nonumber\\[5pt]
&=\frac{2-2x^{m+2}y^{m+2}-(x+y)(1-x^{m+1}y^{m+1})-(1-xy)(x^{m+1}+y^{m+1})}
{(1-xy)(1-x)^{m+1}(1-y)^{m+1}} \nonumber\\[5pt]
&\quad -\frac{1-x^{m+1}y^{m+1}}{(1-xy)(1-x)^{m}(1-y)^{m}}.
\end{align*}
After simplification, we obtain \eqref{eq:1x1yrs}. \qed

\medskip

It is easy to see that \eqref{eq:1x1yrs} may be written as:
\begin{align}
&\hskip -3mm (1-x^m)(1-y^m)\nonumber\\[5pt]
&=\sum_{k=0}^{m-1}\sum_{i=0}^{k}
\frac{m}{m-k}{m-k+i-1\choose i}{m-i-1\choose k-i}x^iy^{k-i}(1-x)^{m-k}(1-y)^{m-k}.
\label{eq:1x1yik}
\end{align}

\noindent{\it Remark.} Applying the multivariate Lagrange inversion
formula, we can also prove \eqref{eq:mrsxy0} and \eqref{eq:1x1yrs}
as well as the following generalization of \eqref{eq:mrsxy}:
\begin{align*}
\sum_{r_1,\ldots,r_m\leq n}
\prod_{k=1}^{m}{n-r_k\choose r_{k+1}}
\frac{x_k^{r_k}}
{(1-x_k)^{r_k+r_{k+1}}}
=\frac{1-(-1)^{m(n+1)}x_1^{n+1}\cdots x_m^{n+1}}
{1-(-1)^{m}x_1\cdots x_m}\prod_{k=1}^{m}\frac{1}{(1-x_k)^n},
\end{align*}
where $r_{m+1}=r_1$.

\medskip

Recall the Vandermonde determinant formula:
\begin{align}
\det(x_i^{n-j})_{1\leq i,j\leq n}=\prod_{1\leq i<j\leq n}(x_i-x_j).
\end{align}
Let $e_i(x_1,\ldots,x_n)$ ($0\leq i\leq n$) be the $i$-th
\emph{elementary symmetric function} of
$x_1,\ldots,x_n$, and  let
$$
(x_1,\ldots,\hat{x_j},\ldots, x_n)=(x_1,\ldots,x_{j-1}, x_{j+1},\ldots, x_n),\qquad
1\leq j\leq n.
$$
\begin{lem}\label{lem:ainv}
Let $A=(x_i^{n-j})_{1\leq i,j\leq n}$ be  the Vandermonde matrix. Then
$$
A^{-1}=\left((-1)^{n-i}\frac{e_{i-1}(x_1,\ldots,\hat{x_j},\ldots, x_n)}
{\prod_{k=1,k\neq j}^n(x_k-x_j)}\right)_{1\leq i,j\leq n}.
$$
\end{lem}
\pf The elementary symmetric functions satisfy the identity
$$
\sum_{k=0}^{n}(-t)^{n-k} e_k(x_1,\ldots,x_n) =\prod_{k=1}^{n}(x_k-t).
$$
Therefore, for each $j=1,2,\ldots, n$, we have
\begin{align}
\sum_{k=1}^{n}(-t)^{n-k} e_{k-1}(x_1,\ldots,\hat{x_j},\ldots, x_n)
=\prod_{\substack{k=1\\ k\neq j}}^{n}(x_k-t). \label{eq:ekxij}
\end{align}
The result then follows by  setting $t=x_{i}$ ($1\leq i\leq n$)
in \eqref{eq:ekxij}.  \qed

\medskip

We shall need the following variant of Vandermonde's determinant.
\begin{lem}\label{lem:xi}
Let $A=((1-x^i)^{n+1-j}(1-x^{i+1})^{n+1-j}x^{i(j-1)})_{1\leq i,j\leq n}$.
Then
\begin{align*}
\det A=(-1)^{\frac{n(n-1)}{2}}x^{\frac{n(n^2-1)}{6}}
\prod_{k=1}^{n}(1-x^{2k-1})^{n+1-k} (1-x^{2k})^{n+1-k}.
\end{align*}
\end{lem}
\pf Extracting $x^{(n-1)i}(1-x^i)(1-x^{i+1})$ from the $i$-th row ($1\leq i\leq n$)
of $A$  and then applying the Vandermonde determinant formula, we obtain
\begin{align*}
\det A&=\prod_{i=1}^nx^{(n-1)i}(1-x^i)(1-x^{i+1})
\cdot\det\left(\frac{(1-x^i)^{n-j}(1-x^{i+1})^{n-j}}
{x^{i(n-j)}}\right)_{1\leq i,j\leq n} \\[5pt]
&=\prod_{i=1}^{n}x^{(n-1)i}(1-x^i)(1-x^{i+1})
\cdot\prod_{1\leq i< j\leq n}
\left(\frac{(1-x^i)(1-x^{i+1})}{x^i}-\frac{(1-x^j)
(1-x^{j+1})}{x^j}\right)\\[5pt]
&=\prod_{i=1}^{n}x^{(n-1)i}(1-x^i)(1-x^{i+1})
\cdot\prod_{1\leq i< j\leq n}
\frac{-(1-x^{j-i})(1-x^{i+j+1})}{x^{j}},
\end{align*}
which yields the desired formula after simplification. \qed

\section{Proof of Theorems}\label{appl}
\begin{thm}\label{thm:pmky}
For any $m\in \mathbb{N}$, there exist polynomials
$P_{m,k}(y)\in\mathbb{Z}[y]$ such that
\begin{align}
&\hskip -3mm \sum_{k=0}^{m}(-1)^k\left[{2m\choose k}-{2m\choose k-2}\right]
\frac{(1-x^{m+1-k})(1-x^{m+1-k}y^{m+1-k})x^k}{1-y^{m+1-k}} \nonumber\\[5pt]
&=\sum_{k=0}^{m}(-1)^k P_{m,k}(y)\frac{(1-x)^{m+1-k}(1-xy)^{m+1-k}(1-y)^{3k}x^k}
{\prod_{i=0}^{k}(1-y^{m+1-i})}.  \label{eq:p2mk}
\end{align}
\end{thm}

\pf
By formula \eqref{eq:1x1yik}, we have
\begin{align*}
&\hskip -3mm 
{(1-x^{m+1-k})(1-x^{m+1-k}y^{m+1-k})x^k} \\[5pt]
&=\sum_{r=0}^{m-k}\sum_{i=0}^{r}\frac{m-k+1}{m-k-r+1}
  {m-k-r+i\choose i}{m-k-i\choose r-i} \\[5pt]
&\quad{}\times {x^{k+r} y^{r-i}(1-x)^{m-k-r+1}(1-xy)^{m-k-r+1}}.
\end{align*}
Therefore, setting $s=r+k$, we obtain
\begin{align}
&\hskip -3mm\sum_{k=0}^{m}\frac{(-1)^k}{1-y^{m+1-k}}
\left[{2m\choose k}-{2m\choose k-2}\right]
(1-x^{m+1-k})(1-x^{m+1-k}y^{m+1-k})x^k\nonumber\\[5pt]
&=\sum_{s=0}^{m}\overline{P}_{m,s}(y)\frac{(1-x)^{m+1-s}(1-xy)^{m+1-s}x^s}
{\prod_{i=0}^{k}(1-y^{m+1-i})},\label{eq:pbar}
\end{align}
where $\overline{P}_{m,s}(y)$ are polynomials given by
\begin{align}
\overline{P}_{m,s}(y)&=\sum_{k=0}^{s}(-1)^k\left[{2m\choose k }-{2m\choose k -2}\right]
\prod_{\substack{j=0\\ j\neq k}}^{s}{(1-y^{m+1-j})}
\nonumber\\[5pt]
&\quad{}\times \sum_{i=0}^{s-k}\frac{m-k+1}{m-s+1}
{m-s+i\choose i}{m-k-i\choose s-k-i}y^{s-k-i}.
\label{eq:pbarms}
\end{align}
By \eqref{eq:ppos}, we have $\overline{P}_{m,s}(y)\in\mathbb{Z}[y]$.
It remains to show that $(1-y)^{3s}\,|\,\overline{P}_{m,s}(y)$.

In view of Lemma \ref{lem:main}, setting $x=q^n$ and $y=q$ in \eqref{eq:pbar},
the left-hand side reduces to $(1-q^2)(1-q)^{2m} S_{2m+1,n}(q)$.
Therefore, it follows from \eqref{eq:qpower} and
\eqref{eq:pbar} that
\begin{align}
\sum_{k=0}^{m}\overline{P}_{m,k}(q)\frac{(1-q^n)^{m+1-k}(1-q^{n+1})^{m+1-k}q^{kn}}
{\prod_{i=0}^{k}(1-q^{m+1-i})}
=\sum_{k=1}^{n}(1-q^{2k})(1-q^k)^{2m}q^{(m+1)(n-k)}. \label{eq:xy2q}
\end{align}
Taking $n=1,2,\ldots,m+1$ in \eqref{eq:xy2q}, we obtain the following matrix
equation
\begin{align}\label{eq:axb}
A\cdot \begin{pmatrix}x_1\\ x_2 \\ \vdots \\ x_{m+1}  \end{pmatrix}
=\begin{pmatrix}b_1\\ b_2 \\ \vdots \\ b_{m+1}  \end{pmatrix},
\end{align}
where $A=((1-q^i)(1-q^{i+1})q^{mi}a_i^{m+1-j})_{1\leq i,j\leq m+1}$
with $a_i=(1-q^i)(1-q^{i+1})/q^i$, and where
\begin{align}
b_i&=\sum\limits_{k=1}^{i}(1-q^{2k})(1-q^k)^{2m}q^{(m+1)(i-k)},\nonumber\\[10pt]
x_j&=\displaystyle\frac{\overline{P}_{m,j-1}(q)}{\prod_{k=0}^{j-1}(1-q^{m+1-k})}.
\label{eq:*}
\end{align}
Now Lemma \ref{lem:xi} implies that $\det A\neq 0$, so Equation \eqref{eq:axb} has
a unique solution given by
\begin{equation}\label{eq:sys}
x_j=\sum_{i=1}^{m+1}(A^{-1})_{ji}b_i, \quad j=1,2,\ldots, m+1,
\end{equation}
where, by Lemma \ref{lem:ainv},
\begin{align*}
(A^{-1})_{ji}&=\frac{(-1)^{m+1-j}}{q^{mi}(1-q^i)(1-q^{i+1})}
\frac{e_{j-1}(a_1,\ldots,\hat{a_i},\ldots, a_{m+1})}
{\prod_{k=1,k\neq i}^{m+1}(a_k-a_i)}\\[5pt]
&=\frac{(-1)^{m+i-j}
q^{m+2-i\choose 2}(1-q^{2i+1})e_{j-1}(a_1,\ldots,\hat{a_i},\ldots, a_{m+1})}
{(q;q)_{m+i+2}(q;q)_{m-i+1}}.
\end{align*}
Here we have adopted the notation $(q;q)_n=(1-q)(1-q^2)\cdots (1-q^n)$.
It follows from \eqref{eq:*} and \eqref{eq:sys} that
\begin{align*}
\overline{P}_{m,j-1}(q)
=\sum_{i=1}^{m+1}(-1)^{m+i-j}q^{m+2-i\choose 2}L_{ij}(q),
\end{align*}
where 
\begin{align}
L_{ij}(q)
=\frac{(1-q^{2i+1})
e_{j-1}(a_1,\ldots,\hat{a_i},\ldots, a_{m+1}) (q;q)_{m+1}b_i}
{(q;q)_{m+i+2}(q;q)_{m-i+1}(q;q)_{m-j+1}}.
\label{eq:pmjba}
\end{align}
Since $a_i={(1-q^i)(1-q^{i+1})}/{q^i}$,
 the valuation of $(1-q)$ in $e_{j-1}(a_1,\ldots,\hat{a_i},\ldots, a_{m+1})$
is at least $2j-2$. Also, it is clear that $(1-q)^{2m+1}\,|\,b_i$ 
for $i=1,2,\ldots,m+1$. Hence, from \eqref{eq:pmjba}
it follows that the valuation of $(1-q)$ in $L_{ij}(q)$
is at least $3j-3$ for $j=1,2,\ldots, m+1$.
Therefore, the polynomials
\begin{align}
{P}_{m,j}(q)=(-1)^j\frac{\overline{P}_{m,j}(q)}{(1-q)^{3j}},\quad j=0,1,\ldots,m,
\label{eq:pmj-1}
\end{align}
satisfy \eqref{eq:p2mk}. This completes the proof. \qed

\medskip

The proof of Theorem \ref{thm:pmn} then follows from \eqref{eq:xy2q} and \eqref{eq:pmj-1}.

Multiplying \eqref{eq:mrsxy} by $(1-x)^m (1-y)^m$ and putting $r=i$ and $s=k-i$, we get
\begin{align*}
\frac{1-x^{m+1}y^{m+1}}{1-xy}=\sum_{k=0}^{m}\sum_{i=0}^{k} {m-k+i\choose i}{m-i\choose k-i}
x^iy^{k-i}(1-x)^{m-k}(1-y)^{m-k}.
\end{align*}
Writing
\begin{align*}
\frac{1-x^{2m+1}y^{2m+1}}{1-xy}
=\frac{1-x^{2m+2}y^{2m+2}}{1-x^2y^2}+xy\frac{1-x^{2m}y^{2m}}{1-x^2y^2},
\end{align*}
we obtain
\begin{align}
&\hskip -3mm \frac{1-x^{2m+1}y^{2m+1}}{1-xy}\nonumber\\[5pt]
&=\sum_{k=0}^{m}x^k\sum_{i=0}^{k}
{m-k+i\choose i}{m-i\choose k-i}y^{2k-2i}(1-x)^{m-k}(1-xy^2)^{m-k}\nonumber\\[5pt]
&\quad +\sum_{k=1}^{m}x^k\sum_{i=0}^{k-1}
{m-k+i\choose i}{m-i-1\choose k-i-1}y^{2k-2i-1}(1-x)^{m-k}(1-xy^2)^{m-k}. \label{eq:xy2m1}
\end{align}

\begin{thm}
For any $m\in \mathbb{N}$, there exist polynomials
$Q_{m,k}(y)\in\mathbb{Z}[y]$ such that
\begin{align}
&\hskip -3mm \sum_{k=0}^{m}(-1)^k\left[{2m-1\choose k}-{2m-1\choose k-2}\right]
\frac{(1-x^{2m-2k+1}y^{2m-2k+1})x^{k}}{1-y^{2m-2k+1}} \nonumber\\[5pt]
&=\sum_{k=0}^{m}(-1)^k Q_{m,k}(y)\frac{(1-x)^{m-k}(1-xy^2)^{m-k}(1-xy)(1-y)^{3k}x^{k}}
{\prod_{i=0}^{k}(1-y^{2m-2i+1})}. \label{eq:qmj-1}
\end{align}
\end{thm}

\pf
By formula \eqref{eq:xy2m1}, we have
\begin{align*}
&\hskip -3mm 
{(1-x^{2m-2k+1}y^{2m-2k+1})x^{k}}\\[5pt]
&= 
\sum_{r=0}^{m-k}\sum_{i=0}^{r}
{m-k-r+i\choose i}\bigg[{m-k-i\choose r-i}y^{2r-2i} 
+{m-k-i-1\choose r-i-1}y^{2r-2i-1}\bigg] \\[5pt]
&\quad{}\times {x^{k+r} (1-x)^{m-k-r}(1-xy^2)^{m-k-r}(1-xy)}.
\end{align*}
Therefore, setting $s=r+k$, we get
\begin{align*}
&\hskip -3mm
\sum_{k=0}^{m}(-1)^k\left[{2m-1\choose k}-{2m-1\choose k-2}\right]
\frac{(1-x^{2m-2k+1}y^{2m-2k+1})x^{k}}{1-y^{2m-2k+1}}  \\[5pt]
&=\sum_{s=0}^{m}\overline{Q}_{m,s}(y)
\frac{x^{s} (1-x)^{m-s}(1-xy^2)^{m-s}(1-xy)}{\prod_{i=1}^{k}(1-y^{2m-2i+1})},
\end{align*}
where
\begin{align}
\overline{Q}_{m,s}(y)
&=\sum_{k=0}^{s}(-1)^k\left[{2m-1\choose k}-{2m-1\choose k-2}\right]
\prod_{\substack{j=0\\ j\neq k}}^{s}(1-y^{2m-2j+1})
\sum_{i=0}^{s-k} {m-s+i\choose i}\nonumber\\[5pt]
&\quad{}\times\bigg[{m-k-i\choose s-k-i}y^{2s-2k-2i} +{m-k-i-1\choose s-k-i-1}y^{2s-2k-2i-1}\bigg]. \label{eq:qbarms}
\end{align}
What remains is to show that
$$(1-q)^{s}(1-q^2)^{2s}\,|\,\overline{Q}_{m,s}(q) \text{ and }
{Q}_{m,s}(q)=(-1)^s\frac{\overline{Q}_{m,s}(q)}{(1-q)^{s}(1-q^2)^{2s}}.$$
The proof is exactly the same as that of Theorem \ref{thm:pmky} and is omitted. \qed

\medskip

The proof of Theorem \ref{thm:qmn} then follows
from \eqref{eq:oddq} and \eqref{eq:qmj-1}
(replacing $x$ and $y$ by $q^n$ and $q^{\frac{1}{2}}$, respectively).

\medskip
\noindent{\it Remark.} It follows from \eqref{eq:cor1} and \eqref{eq:cor2}
that $P_{m,m}(q)=Q_{m,m}(q)=0$
if $m\geq 1$, and $P_{m,s}(0)=Q_{m,s}(0)={m+s-2\choose s}-{m+s-2\choose s-2}$.
Moreover, if $s<m$,   the polynomial $P_{m,s}(q)$ has degree $s(2m-3-s)/2$
while $Q_{m,s}(q)$ has degree $s(2m-3-s)$.

\medskip
\noindent{\it Proof of Theorem \ref{thm:salie}.}
By formula \eqref{eq:1x1yik}, we have
\begin{align*}
{(1-x^{m-r})(1-x^{m-r}y^{m-r})x^r}
&=\sum_{k=0}^{m-r}\sum_{i=0}^{k}\frac{m-r}{m-r-k} 
{m-r-k+i-1\choose i}{m-r-i-1\choose k-i}
\\[5pt]
&\quad{}\times 
{x^{r+k}y^{k-i}(1-x)^{m-r-k}(1-xy)^{m-r-k}}.
\end{align*}
Therefore, setting $s=r+k$ and
\begin{align*}
&\hskip -2mm \overline{G}_{m,s}(y) \\[5pt]
&=\sum_{r=0}^{s}(-1)^r {2m\choose r}
\prod_{\substack{j=0\\ j\neq r}}^{s}{(1+y^{m-j})}
\sum_{i=0}^{s-r}\frac{m-r}{m-s}
{m-s+i-1\choose i}{m-r-i-1\choose s-r-i}y^{s-r-i},
\end{align*}
we obtain
\begin{align}
\sum_{r=0}^{m-1}(-1)^r{2m\choose r}
\frac{(1-x^{m-r})(1-x^{m-r}y^{m-r})x^r}{1+y^{m-r}}
=\sum_{s=0}^{m}\overline{G}_{m,s}(y)\frac{(1-x)^{m-s}(1-xy)^{m-s}x^s}
{\prod_{i=0}^{s}(1+y^{m-i})}.\label{eq:gmsbar}
\end{align}
Similarly to the proof of Theorem \ref{thm:pmky}, we can show that
$G_{m,s}(y)=(-1)^s \overline{G}_{m,s}(y)/(1-y)^{2s}$ is a polynomial in $\mathbb{Z}[y]$.

Theorem~\ref{thm:salie} then follows from Lemma~\ref{lem:sum} after
substituting $x=q^n$ and $y=q$ into \eqref{eq:gmsbar}.  \qed

\medskip
The proof of Theorem~\ref{thm:salodd} is analogous to that of Theorem \ref{thm:qmn}
and is omitted here.

\section{Sums of $m$-th Powers for $m\leq 11$}\label{compute}
Theorems \ref{thm:pmn}--\ref{thm:salodd} permit us to
compute $P_{m,k}(q)$, $Q_{m,k}(q)$, $G_{m,k}(q)$, and $H_{m,k}(q)$ quickly by using Maple.
Tables \ref{t:p}--\ref{t:hmk} give the first values of these polynomials.

\begin{table}[h]
\caption{Values of $P_{m,k}(q)$ for $0\leq m\leq 5$.\label{t:p}}
{\footnotesize
\begin{center}
\begin{tabular}{|l|r|r|r|r|r|r|}
\hline
$k\setminus m$ &0 & 1 & 2 & 3 & 4 & 5\\\hline
0              & 1  & 1 & 1 & 1 & 1 & 1\\\hline
1              &&   & 1 & $2(q+1)$ & $3q^2+4q+3$ & $2(q+1)(2q^2+q+2)$\\\hline
2              &&   &   & $2(q+1)$ & $(q+1)(5q^2+8q+5)$ & $(q+1)(9q^4+19q^3+29q^2+19q+9)$\\\hline
3              &&   &   &          & $(q+1)(5q^2+8q+5)$ & $2(q+1)^2(q^2+q+1)(7q^2+11q+7)$\\\hline
4              &&   &   &          &                    & $2(q+1)^2(q^2+q+1)(7q^2+11q+7)$\\\hline
\end{tabular}
\end{center}
}
\end{table}
\vskip -5mm

\begin{table}[h]
\caption{Values of $Q_{m,k}(q)$ for $1\leq m\leq 4$.\label{t:q}}
{\footnotesize
\begin{center}
\begin{tabular}{|l|r|r|r|r|}
\hline
$k\setminus m$ & 1 & 2 & 3 & 4 \\\hline
0              & 1 & 1 & 1 & 1 \\\hline
1              &   & 1     &  $2q^2+q+2$  & $3q^4+2q^3+4q^2+2q+3$
\\\hline
2              &   &   & $2q^2+q+2$ & $(q^2+q+1)(5q^4+q^3+9q^2+q+5)$
\\\hline
3              &   &   &          & $(q^2+q+1)(5q^4+q^3+9q^2+q+5)$ \\\hline
\end{tabular}
\end{center}
}
\end{table}
For $m=5$, we have $Q_{5,0}(q)=1$, and
{\footnotesize 
\begin{align*}
Q_{5,1}(q)&=4q^6+3q^5+6q^4+4q^3+6q^2+3q+4, \\
Q_{5,2}(q)&=9q^{10}+13q^{9}+33q^{8}+37q^{7}+61q^{6}+51q^{5}+61q^{4}+37q^{3}+33q^{2}+13q+9,
\\
Q_{5,3}(q)&=Q_{5,4}(q) \\
&=(q^2+q+1)(14q^{10}+14q^{9}+56q^{8}+46q^{7}+100q^{6}+65q^{5}
+100q^{4}+46q^{3}+56q^{2}+14q+14).
\end{align*}
}
\vskip -5mm

\begin{table}[h]
\caption{Values of $G_{m,k}(q)$ for $1\leq m\leq 5$.\label{t:gmk}}
{\footnotesize
\begin{center}
\begin{tabular}{|l|r|r|r|r|r|}
\hline
$k\setminus m$ & 1 & 2 & 3 & 4 & 5\\\hline
0              & 1 & 1 & 1 & 1 & 1\\\hline
1              &   & 2 & $3(q+1)$&$4(q^2+q+1)$ & $5(q+1)(q^2+1)$\\\hline
2              &   &   & $6(q+1)$&$2(q+1)(5q^2+7q+5)$&$5(q+1)(3q^4+4q^3+8q^2+4q+3)$\\\hline
3              &   &   &         &$4(q+1)(5q^2+7q+5)$&$5(q+1)^2(7q^4+14q^3+20q^2+14q+7)$\\\hline
4              &   &   &         &                   &$10(q+1)^2(7q^4+14q^3+20q^2+14q+7)$\\\hline
\end{tabular}
\end{center}
}
\end{table}

\vskip -5mm
\begin{table}[!h]
\caption{Values of $H_{m,k}(q)$ for $1\leq m\leq 4$.\label{t:hmk}}
{\footnotesize
\begin{center}
\begin{tabular}{|l|r|r|r|r|}
\hline
$k\setminus m$&1 & 2 &     3          & 4   \\\hline
0             &1 & 1 &     1          & 1  \\\hline
1             &  & 2 & $3q^2+2q+3$    & $4q^4+3q^3+4q^2+3q+4$ \\\hline
2             &  &   & $2(3q^2+2q+3)$ & $10q^6+15q^5+30q^4+26q^3+30q^2+15q+10$ \\\hline
3             &  &   &                & $2(10q^6+15q^5+30q^4+26q^3+30q^2+15q+10)$  \\\hline
\end{tabular}
\end{center}
}
\end{table}

Substituting the values of Tables \ref{t:p} and \ref{t:q} into
Theorems \ref{thm:pmn} and \ref{thm:qmn} yields the
summation formulas for sums of $m$-th power for $m=1,2,\ldots, 11$.
In particular, for $1\leq m\leq 5$ we recover the formulas 
\eqref{eq:f1}--\eqref{eq:f5}
of Warnaar and Schlosser. For $m=6,7,\ldots, 11$ we obtain
the following formulas of Faulhaber type:
\begin{align*}
S_{6,n}(q)&=\frac{(1-q^n)(1-q^{n+1})(1-q^{n+\frac{1}{2}})}
{(1-q)(1-q^2)(1-q^{\frac{7}{2}})}
\Bigg[\frac{(1-q^n)^2(1-q^{n+1})^2}{(1-q)^4}
\\[5pt]
&\quad{}-(2+2q+q^{\frac{1}{2}})\left(\frac{(1-q^n)(1-q^{n+1})q^{n}}
{(1+q^{\frac{1}{2}})(1-q)(1-q^{\frac{5}{2}})}
-\frac{(1-q^{\frac{1}{2}})^2q^{2n}}
{(1-q^{\frac{3}{2}})(1-q^{\frac{5}{2}})}\right)\Bigg],
\\[10pt]
S_{7,n}(q)&=\frac{(1-q^n)^2(1-q^{n+1})^2}{(1-q)^3(1-q^4)}
\\[8pt]
&\quad{}\times\Bigg[\frac{(1-q^n)^2(1-q^{n+1})^2}{(1-q)^3(1-q^2)}
-\frac{2(1-q^n)(1-q^{n+1})q^{n}}{(1-q)(1-q^3)}
+\frac{2(1-q)^2 q^{2n}}{(1-q^2)(1-q^3)}\Bigg],
\end{align*}
\begin{align*}
S_{8,n}(q)&=\frac{(1-q^n)(1-q^{n+1})(1-q^{n+\frac{1}{2}})}
{(1-q)(1-q^2)(1-q^{\frac{9}{2}})}
\\[5pt]
&\quad{}\times \Bigg[\frac{(1-q^n)^3(1-q^{n+1})^3}{(1-q)^6}
-(3+2q^{\frac{1}{2}}+4q+2q^{\frac{3}{2}}+3q^2)
\frac{(1-q^n)^2(1-q^{n+1})^2q^{n}}
{(1+q^{\frac{1}{2}})(1-q)^3(1-q^{\frac{7}{2}})}
\\[5pt]
&\quad{}+(5+q^{\frac{1}{2}}+9q+q^{\frac{3}{2}}+5q^2)
\frac{(1-q^{\frac{3}{2}})(1-q^n)(1-q^{n+1})q^{2n}}
{(1+q^{\frac{1}{2}})(1-q)(1-q^{\frac{5}{2}})(1-q^{\frac{7}{2}})}\\[5pt]
&\quad{}-(5+q^{\frac{1}{2}}+9q+q^{\frac{3}{2}}+5q^2)
\frac{(1-q^{\frac{1}{2}})^2q^{3n}}{(1-q^{\frac{5}{2}})(1-q^{\frac{7}{2}})}\Bigg],
\\[10pt]
S_{9,n}(q)&=\frac{(1-q^n)^2(1-q^{n+1})^2}{(1-q)^3(1-q^5)}\\[5pt]
&\quad{}\times
\Bigg[\frac{(1-q^n)^3(1-q^{n+1})^3}{(1-q)^5(1-q^2)}
-\frac{(3q^2+4q+3)(1-q^n)^2(1-q^{n+1})^2 q^{n}}{(1-q)^2(1-q^2)(1-q^4)}
\\[5pt]
&\quad{}+(5q^2+8q+5)\left(\frac{(1-q^n)(1-q^{n+1}) q^{2n}}{(1-q^3)(1-q^4)}
-\frac{(1-q)^3 q^{3n}}{(1-q^2)(1-q^3)(1-q^4)}\right)\Bigg], 
\end{align*}
\begin{align*}
&S_{10,n}(q)\\[5pt]
&=\frac{(1-q^n)(1-q^{n+1})(1-q^{n+\frac{1}{2}})}
{(1-q)(1-q^2)(1-q^{\frac{11}{2}})}
\Bigg[\frac{(1-q^n)^4(1-q^{n+1})^4}{(1-q)^8}
\\[5pt]
&\quad{}-Q_{5,1}(q^{\frac1 2})
 \frac{(1-q^{\frac1 2})(1-q^n)^3(1-q^{n+1})^3 q^{n}}{(1-q)^6(1-q^{\frac{9}{2}})}
+Q_{5,2}(q^{\frac1 2})\frac{(1-q^{\frac1 2})^2(1-q^n)^2(1-q^{n+1})^2q^{2n}}
{(1-q)^4(1-q^{\frac{7}{2}})(1-q^{\frac{9}{2}})}\\[5pt]
&\quad{}-Q_{5,3}(q^{\frac1 2})
 \frac{(1-q^{\frac1 2})^3(1-q^n)(1-q^{n+1})q^{3n}}
{(1-q)^2(1-q^{\frac{5}{2}})(1-q^{\frac{7}{2}})(1-q^{\frac{9}{2}})}
+\frac{Q_{5,4}(q^{\frac1 2})(1-q^{\frac1 2})^4 q^{4n}}
{(1-q^{\frac{3}{2}})(1-q^{\frac{5}{2}})(1-q^{\frac{7}{2}})(1-q^{\frac{9}{2}})}
\Bigg],\\[10pt]
&S_{11,n}(q)\\[5pt]
&=\frac{(1-q^n)^2(1-q^{n+1})^2}{(1-q)^3(1-q^6)}
\Bigg[\frac{(1-q^n)^4(1-q^{n+1})^4}{(1-q)^7(1-q^2)}
-\frac{2(2q^2+q+2)(1-q^n)^3(1-q^{n+1})^3 q^{n}}{(1-q)^5(1-q^5)}
\\[5pt]
&\quad{}+
\frac{(9q^4+19q^3+29q^2+19q+9)(1-q^n)^2(1-q^{n+1})^2 q^{2n}}{(1-q)^2(1-q^4)(1-q^5)}
\\[5pt]
&\quad{}-2(q+1)(7q^2+11q+7)\left(\frac{(1-q^n)(1-q^{n+1})q^{3n}}{(1-q^4)(1-q^5)}
-\frac{(1-q)^3q^{4n}}{(1-q^2)(1-q^4)(1-q^5)}\right)\Bigg].
\end{align*}
{}From the computational point of view, with the help of Maple or
other softwares, it is, of course, not difficult to give further extension of
the above list of $S_{m,n}(q)$'s.

\section{Further Remarks}\label{sec:open}
For $r\in \mathbb{N}$, define the following more general summation
$$
S_{m,n,r}(q)=\sum_{k=1}^{n}\frac{1-q^{(2r+2)k}}{1-q^{2r+2}}
\left(\frac{1-q^k}{1-q}\right)^{m-1}q^{\frac{m+2r+1}{2}(n-k)}.
$$
Then we can also obtain a similar summation formula.
\begin{thm}\label{thm:pr}
For $0\leq r\leq m$, there exist polynomials
$P_{m,k,r}(q)\in\mathbb{Z}[q]$ such that
\begin{align*}
S_{2m-2r+1,n,r}(q) &=\sum_{k=0}^{m}(-1)^kP_{m,k,r}(q)
\frac{(1-q^n)^{m+1-k}(1-q^{n+1})^{m+1-k}q^{kn}}
{(1-q^{2r+2})(1-q)^{2m-3k}\prod_{i=0}^{k}(1-q^{m+1-i})}.
\end{align*}
Furthermore, we have
\begin{align*}
P_{m,s,r}(q)
&=\frac{\prod_{j=0}^{s}(1-q^{m+1-j})}{(1-q)^{3s-2r}}
\sum_{k=0}^{s}\frac{(-1)^{s-k}}{1-q^{m+1-k}}
\left[{2m-2r\choose k }-{2m-2r\choose k-2r-2}\right]
  \\[5pt]
&\quad{}\times\sum_{i=0}^{s-k}\frac{m-k+1}{m-s+1}
{m-s+i\choose i}{m-k-i\choose s-k-i} q^{s-k-i}.
\end{align*}
\end{thm}
\begin{thm}\label{thm:qr}
For $0\leq r\leq m$, there exist polynomials
$Q_{m,k,r}(q)\in\mathbb{Z}[q]$ such that
\begin{align*}
S_{2m-2r,n,r}(q)
=\sum_{k=0}^{m}(-1)^k Q_{m,k,r}(q^{\frac{1}{2}})
\frac{(1-q^{n+\frac{1}{2}})(1-q^n)^{m-k}(1-q^{n+1})^{m-k}
(1-q^{\frac{1}{2}})^k q^{kn}}
{(1-q^{2r+2})(1-q)^{2m-2k-1}\prod_{i=0}^{k}(1-q^{m-i+\frac{1}{2}})}.
\end{align*}
Furthermore, we have
\begin{align*}
Q_{m,s,r}(q) &=\frac{\prod_{j=0}^{s}(1-q^{2m-2j+1})}{(1-q)^s(1-q^2)^{2s-2r}} 
\sum_{k=0}^{s}\frac{(-1)^{s-k}}{1-q^{2m-2k+1}}
\left[{2m-2r-1\choose k}-{2m-2r-1\choose k-2r-2}\right]
 \\[5pt]
&\quad{}\times\sum_{i=0}^{s-k} {m-s+i\choose i}\bigg[{m-k-i\choose s-k-i}q^{2s-2k-2i}
 +{m-k-i-1\choose s-k-i-1}q^{2s-2k-2i-1}\bigg].
\end{align*}
\end{thm}

For example, we have
\begin{align*}
S_{4,n,1}(q)&=\frac{(1-q^n)^3(1-q^{n+1})^3(1-q^{n+\frac{1}{2}})}
{(1-q)^3(1-q^4)(1-q^{\frac{7}{2}})}
+\frac{(1+4q^{\frac{1}{2}}+4q+4q^{\frac{3}{2}}+q^2)(1-q^{\frac{1}{2}})}{1-q}\\[5pt]
&\quad{}\times\frac{(1-q^n)(1-q^{n+1})(1-q^{n+\frac{1}{2}})}
{(1-q^{\frac{5}{2}})(1-q^{\frac{7}{2}})(1-q^4)}
\Bigg[\frac{(1-q^n)(1-q^{n+1})}{(1-q)^2}q^{n+\frac{1}{2}}
-\frac{1-q^{\frac{1}{2}}}{1-q^{\frac{3}{2}}}
q^{2n+\frac{1}{2}}\Bigg], 
\end{align*}
and
\begin{align*}
S_{5,n,1}(q)&=\frac{(1-q^n)^4(1-q^{n+1})^4}{(1-q)^4(1-q^4)^2}\\[5pt]
&\quad{}+\frac{4(1-q^n)^2(1-q^{n+1})^2}{(1-q)(1-q^3)(1-q^4)^2}
\Bigg[\frac{(1+q)(1-q^n)(1-q^{n+1})}{(1-q)^2}q^{n+1}-q^{2n+1}\Bigg].
\end{align*}

There is a similar formula for
$$
T_{m,n,r}(q)=\sum_{k=1}^{n}(-1)^{n-k}\frac{1-q^{(2r+1)k}}{1-q^{2r+1}}
\left(\frac{1-q^k}{1-q}\right)^{m-1}q^{\frac{m+2r}{2}(n-k)},
$$
which is left to the interested readers.

In a forthcoming paper \cite{GRZ}, it will be shown that the
coefficients of the polynomials $P_{m,k}(q)$, $Q_{m,k}(q)$,
$G_{m,k}(q)$ and $H_{m,k}(q)$ are actually \emph{nonnegative integers}
and have interesting combinatorial interpretations in terms
of nonintersecting lattice paths.

\vskip 5mm
\noindent{\bf Acknowledgment.}
The second author was supported by EC's IHRP Programme, within Research Training
Network ``Algebraic Combinatorics in Europe," grant HPRN-CT-2001-00272.

\renewcommand{\baselinestretch}{1}

\end{document}